\documentstyle[12pt]{article}
\begin{document}
\vskip 1.5in
\centerline{\bf Some Integrals Involving Bessel Functions}
\vskip 1.in
\centerline{\bf M.L. Glasser}

Department of Physics and Department of Mathematics and Computer
Science, Clarkson University, Potsdam, N.Y. 13699-5815 (USA)
\vskip .3in
\centerline{\bf E. Montaldi}

Dipartimento di Fisica, Universita di Milano, Via Celoria 16, 20133

Milano, Italy
\vskip 1in

\centerline{\bf Abstract}
\vskip .5in
A number of new definite integrals involving Bessel functions are
presented. These have been derived by finding new integral
representations for the product of two Bessel functions of different
order and argument in terms of the generalized hypergeometric function
with subsequent reduction to special cases. Connection is made with
Weber's second exponential integral and Laplace transforms of products
of three Bessel functions.

\vskip .2in

AMS Classification No. 33A40, 33A30

\vfill\eject
\centerline{\bf Some Integrals Involving Bessel Functions}
\vskip .4in

\centerline{M.L. Glasser* and E. Montaldi**}
\vskip .2in

*Department of Physics and Department of Mathematics and Computer
Science, Clarkson University, Potsdam, N.Y. 13699-5815 (USA)
\vskip .2in

**Dipartimento di Fisica, Universita di Milano, Via Celoria 16, 20133
Milano, Italy

\vskip 1in
{\bf 1. Introduction}

The aim of this work is to derive a number of infinite integrals
involving Bessel functions which appear to be new. The approach is,
beginning with an expression for the product of two Bessel functions as
a sum of Gauss functions, to integrate and perform a resummation to
obtain other hypergeometric functions, and then to reduce these to more
familiar form. In particular we exploit the relation between Bessel
functions and the function $_0F_3$[1]. Finally, we generalize Weber's
second exponential integral by expressing the Laplace transform of a
product of three Bessel functions as an infinite series of products of
modified Bessel functions. We have not aimed at complete rigor or
generality; operations such as  
interchange of limits and Hankel inversion are carried out formally,
i.e. without verifying that the conditions stated ensure their
validity.
In the majority of cases approriate conditions can be supplied by appeal to 
convergence and analytic
continuation.

\vskip .2in
{\bf  2. On the product $J_{\mu}(ax)J_{\nu}(bx)$}

Our starting point is the familiar expansion [2]
$$\Gamma(\mu+1)\Gamma(\nu+1)J_{\mu}(ax)J_{\nu}(bx)=$$
$$(\frac{1}{2}ax)^{\mu}(\frac{1}{2}bx)^{\nu}\sum_{m=0}^{\infty}
\frac{(-1)^m(\frac{1}{2}ax)^{2
m}}{m!(\mu+1)_m}\; _2F_1(-m,-\mu-m
;\nu+1;\frac{b^2}{a^2}) \eqno (2.1)$$
By making use of the standard transformation[3], we have
$$_2F_1(-m,-\mu-m;\nu+1;\frac{b^2}{a^2})=$$$$=(1-\frac{b^2}{a^2})^{\mu+\nu+2m+1}
\;_2F_1(\nu+m+1,\mu+\nu+m+1;\nu+1;\frac{b^2}{a^2}) \eqno (2.2)$$

Next, we recall that [4]
$$_2F_1(\nu + m + 1, \nu+\mu+m+1;\nu+1;\frac{b^2}{a^2})=$$ $$=
\frac{4(a/b)^{\nu}}{\Gamma(\mu+\nu+1)(\nu+1)_m(\mu+\nu+1)_m}\cdot$$
$$\cdot\int_0^{\infty}t^{\mu+\nu+2m+1}K_{\mu}(2t)I_{\nu}(2\frac{b}{a}t)dt \eqno(2.3)$$
$$ m=0,1,2,\dots,Re\;\nu>-1,\;\;\;\;Re\;(\mu+\nu)>-1,\;\;\;\;|Re(\frac{b}{a})|<1.$$
By inserting (2.2) and (2.3) into (2.1), we get
$$\Gamma(\nu+1)\Gamma(\mu+1)\Gamma(\mu+\nu+1)J_{\mu}(ax)J_{\nu}(bx)=
4(\frac{1}{2}ax)^{\mu}(\frac{1}{2}bx)^{\nu}(1-\frac{b^2}{a^2})^{\mu+\nu+1}\cdot $$
$$\cdot(\frac{a}{b})^{\nu}\int_0^{\infty}t^{\mu+\nu+1}\;_0F_3(\mu+1,\nu+1,
\mu+\nu+1;-z^2t^2)K_{\mu}(2t)I_{\nu}(2\frac{b}{a}t)dt \eqno(2.4)$$
$$Re\;\nu>-1,\;\;\;\;Re\;(\mu+\nu)>-1,\;\;\;\;|Re(\frac{b}{a})|<1,$$
where $z=\frac{1}{2}ax(1-\frac{b^2}{a^2})$.

It is interesting to observe that eq. (2.4) enables us to derive the
Mellin transform of $J_{\mu}(ax)J_{\nu}(bx)$, i.e. the
Weber-Schafheitlin integral [5], in a simple way. Indeed from [6],
$$\int_0^{\infty}x^{\mu+\nu-s}\;_0F_3(\mu+1,\nu+1,\mu+\nu+1;-z^2t^2)dx=$$
$$=\frac{1}{2}[\frac{1}{2}a(1-\frac{b^2}{a^2})t]^{s-\mu-\nu-1}
\frac{\Gamma(\frac{\mu+\nu-s+1}{2})\Gamma(\mu+1)\Gamma(\nu+1)\Gamma(\mu+\nu+1)}
{\Gamma(\frac{\mu-\nu+s+1}{2})\Gamma(\frac{\nu-\mu+s+1}{2})\Gamma(\frac{
\mu+\nu+s+1}{2})}, \eqno(2.5)$$
$$Re\;(\mu+\nu-s)>-1,$$
and [4]
\vfill\eject
$$\int_0^{\infty}t^sK_{\mu}(2t)I_{\nu}(2\frac{b}{a}t)dt=\frac{(b/a)^{\nu}}
{4\Gamma(\nu+1)}\Gamma(\frac{\mu+\nu+s+1}{2})\Gamma(\frac{\nu-\mu+s+1}{2})
\cdot$$ $$_2F_1(\frac{\nu+\mu+s+1}{2},\frac{\nu-\mu+s+1}{2};\nu+1;\frac{b^2}{a^2})
 \eqno(2.6)$$
$$Re(\nu\pm\mu+s)>-1,\;\;\;\;|Re(\frac{b}{a})|<1$$
we immediately obtain
$$\int_0^{\infty}x^{-s}J_{\mu}(ax)J_{\nu}(bx)dx=2^{-s}b^{\nu}a^{s-\nu-1}
\frac{\Gamma(\frac{\mu+\nu-s+1}{2})}{\Gamma(\nu+1)\Gamma(\frac{\mu-\nu+s+1}{2})}\cdot
$$
$$(1-\frac{b^2}{a^2})^s\;_2F_1(\frac{\mu+\nu+s+1}{2},
\frac{\nu-\mu+s+1}{2};\nu+1;\frac{b^2}{a^2})$$
$$=2^{-s}b^{\nu}a^{s-\nu-1}\frac{\Gamma(\frac{\mu+\nu-s+1}{2})}{\Gamma(\nu+1)
\Gamma(\frac{\mu-\nu+s+1}{2})}\;_2F_1(\frac{\nu-\mu-s+1}{2},\frac{\nu+\mu-s+1}{2};
\nu+1;\frac{b^2}{a^2})\eqno(2.7)$$
$$Re(\mu+\nu-s)>-1,\;\;\;\; \;\;\;\; 0<  b< a$$

With $b\rightarrow ib\; (b>0)$,and $x$ positive, eq. (2.4) becomes
$$\Gamma(\mu+1)\Gamma(\nu+1)\Gamma(\mu+\nu+1)J_{\mu}(ax)I_{\nu}(bx)=$$
$$(\frac{1}{4}ax)^{\mu}(\frac{1}{4}bx)^{\nu}(1+\frac{b^2}{a^2})^{\mu+\nu+1}
(\frac{a}{b})^{\nu}\cdot$$
$$\int_0^{\infty}t^{\mu+\nu+1}\;_0F_3(\mu+1,\nu+1,\mu+\nu+1;-\frac{1}{16}a^2
x^2(1+\frac{b^2}{a^2})^2t^2)K_{\mu}(t)J_{\nu}(\frac{b}{a}t)dt\eqno(2.8)$$
$$ Re\;\nu>-1,\;\;\;\;Re(\mu+\nu)>-1$$
or, by writing b=ay and x=a/4
$$\Gamma(\mu+1)\Gamma(\nu+1)\Gamma(\mu+\nu+1)J_{\mu}(\frac{1}{4}a^2)
I_{\nu}(\frac{1}{4}a^2y)=(\frac{a^2}{16})^{\mu+\nu}(1+y^2)^{\mu+\nu+1}
\cdot$$
$$\cdot\int_0^{\infty}t^{\mu+\nu+1}\;_0F_3(\mu+1,\nu+1,\mu+\nu+1;-\frac{
a^4}{256}(1+y^2)^2t^2)K_{\mu}(t)J_{\nu}(yt)dt\eqno(2.9)$$
$$Re\;\nu>-1,\;\;\;\;Re(\mu+\nu)>-1.$$
Also, with $a^2\rightarrow\frac{16a}{1+y^2}$,
\vfill\eject
$$\Gamma(\mu+1)\Gamma(\nu+1)\Gamma(\mu+\nu+1)J_{\mu}(\frac{4a}{1+y^2})
I_{\nu}(\frac{4ay}{1+y^2})=(1+y^2)a^{\mu+\nu}\cdot$$
$$\cdot\int_0^{\infty}t^{\mu+\nu+1}\;_0F_3(\mu+1,\nu+1,\mu+\nu+1;-a^2t^2
)K_{\mu}(t)J_{\nu}(yt)dt\eqno(2.10)$$
$$Re\;\nu>-1,\;\;\;\;Re(\mu+\nu)>-1,  a,y>0.$$
For $Re\nu\le -1/2$ we may use the Hankel inversion formula [7], and eq. (2.10) gives
$$(at)^{\mu+\nu}\;_0F_3(\mu+1,\nu+1,\mu+\nu+1;-a^2t^2)K_{\mu}(t)=
\Gamma(\mu+1)\Gamma(\nu+1)\Gamma(\mu+\nu+1)\cdot$$
$$\cdot\int_0^{\infty}\frac{y}{1+y^2}J_{\nu}(ty)J_{\mu}(\frac{4a}{1+y^2})
I_{\nu}(\frac{4ay}{1+y^2})dy \eqno(2.11)$$
This is an interesting addition to the class of Sonine-Gegenbauer
integrals [8]. In particular, by dividing by $a^{\mu+\nu}$ and taking
the limit $a\rightarrow 0$, which may be taken under the integral
sign, we have
$$(\frac{t}{2})^{\mu+\nu}K_{\mu}(t)=\Gamma(\mu+\nu+1)\int_0^{\infty}\frac
{y^{\nu+1}}{(1+y^2)^{\mu+\nu+1}}J_{\nu}(ty)dy \eqno(2.12)$$

The four particular cases
$$\mu=0,\;\;\nu=-\frac{1}{2};\;\;\;\mu=0,\;\;\nu=\frac{1}{2};\;\;\;
\nu=0,\;\;\mu=-\frac{1}{2};\;\;\;\nu=0,\;\;\mu=\frac{1}{2}$$
of eq.(2.9) are of some interest. By using
$$_0F_3(\frac{1}{2},\frac{1}{2},1;-\frac{x^4}{256})=ber(x)\eqno(2.13)$$
and
$$\frac{x^2}{4}\;_0F_3(\frac{3}{2},\frac{3}{2},1:-\frac{x^4}{256})=
bei(x)\eqno(2.14)$$
we obtain the four integral representations
$$\frac{\pi}{2}(1+y^2)^{-1/2}J_0(\frac{a^2}{4})\cosh(\frac{a^2y}{4})=\int_0
^{\infty}K_0(t)ber(a\sqrt{(1+y^2)t})\cos(yt)dt \eqno(2.15)$$
$$\frac{\pi}{2}(1+y^2)^{-1/2}J_0(\frac{a^2}{4})\sinh(\frac{a^2y}{4})=
\int_0^{\infty}K_0(t)bei(a\sqrt{(1+y^2)t})\sin(yt)dt\eqno(2.16)$$
$$(1+y^2)^{-1/2}I_0(\frac{a^2y}{4})\cos(\frac{a^2}{4})=
\int_0^{\infty}e^{-t}ber(a\sqrt{(1+y^2)t})J_0(yt)dt \eqno(2.17)$$
$$(1+y^2)^{-1/2}I_0(\frac{a^2y}{4})\sin(\frac{a^2}{4})=
\int_0^{\infty}e^{-t}bei(a\sqrt{(1+y^2)t})J_0(yt)dt \eqno(2.18).$$
Therefore,
$$K_0(t)ber(a\sqrt{t})=\int_0^{\infty}(1+y^2)^{-1/2}J_0(\frac{1}{4}\frac{
a^2}{1+y^2})\cosh(\frac{1}{4}\frac{a^2y}{1+y^2})\cos(ty)dy\eqno(2.19)
$$
$$K_0(t)bei(a\sqrt{t})=
\int_0^{\infty}(1+y^2)^{-1/2}J_0(\frac{1}{4}\frac{a^2}{1+y^2})\sinh(
\frac{1}{4}\frac{a^2y}{1+y^2})\sin(ty)dy\eqno(2.20)$$
$$\frac{1}{t}e^{-t}ber(a\sqrt{t})=
\int_0^{\infty}y(1+y^2)^{-1/2}I_0(\frac{1}{4}\frac{a^2y}{1+y^2})\cos(\frac
{1}{4}\frac{a^2}{1+y^2})J_0(ty)dy\eqno(2.21)$$
$$\frac{1}{t}e^{-t}bei(a\sqrt{t})=\int_0^{\infty}y(1+y^2)^{-1/2}I_0(
\frac{1}{4}\frac{a^2y}{1+y^2})\sin(\frac{1}{4}\frac{a^2}{1+y^2})
J_0(ty)dy\eqno(2.22)$$
A curious formula arising from (2.4) may be mentioned. We set
$\mu=\frac{1}{2}$, replace $I_{\nu}(2\frac{b}{a}t)$ by Poisson's
integral [9]
$$I_{\nu}(2\frac{b}{a}t)=\frac{1}{\Gamma(\nu+\frac{1}{2})}(\frac
{bt}{a})^{\nu}\int_{-1}^1
e^{-2\frac{b}{a}tu}(1-u^2)^{\nu-\frac{1}{2}}du\eqno(2.23)$$
$$ Re\;\nu>-\frac{1}{2}$$
interchange the order of integration, and then use the Laplace
transform
$$\int_0^{\infty}e^{-\beta t}t^{2\nu+1}\;_0F_3(\frac{3}{2},\nu+1,
\nu+\frac{3}{2};-\alpha t^2)dt=\frac{\Gamma(2\nu+2)}{4\sqrt{\alpha}}
\beta^{-(2\nu+1)}\sin(4\frac{\sqrt{\alpha}}{\beta})\eqno(2.24)$$
$$Re\;\nu>-1,\;\;\;\;Re\;\beta>0$$
which is readily established by expressing the $_0F_3$ as its power
series and integrating term-by-term. This leads to

$$\sin(ax)J_{\nu}(bx)=$$
$$=\frac{1}{\sqrt{\pi}\Gamma(\nu+\frac{1}{2})}
(\frac{1}{2}bx)^{\nu} 
(a^2-b^2)^{\nu+\frac{1}{2}} 
\int_{-1}^1
\frac{(1-u^2)^{\nu-\frac{1}{2}}}{(a+bu)^{2\nu+1}}
\sin(\frac{a^2-b^2}{a+bu}x)du \eqno(2.25)$$
$$Re\;\nu>-\frac{1}{2},\;\;\;\;|a|>|b|>0.$$
Therefore, with $x=\frac{\pi}{2a}$ and $b=ay$
$$J_{\nu}(\frac{\pi y}{2})$$
$$=\frac{1}{\sqrt{\pi}\Gamma(\nu+\frac{1}{2})}(\frac{1}{4}\pi
y)^{\nu}(1-y^2)^{\nu+\frac{1}{2}}\int_{-1}^1
\frac{(1-u^2)^{\nu-1/2}}{(1+uy)^{2\nu+1}}\sin(\frac{\pi}{2}\frac{1-y^2}{
1+uy})du \eqno(2.26)
$$
$$Re\;\nu>-\frac{1}{2},\;\;\;\;|y|<1.$$

By taking advantage of other known reductions for the $_0F_3$ functions
in (2.4)
many additional new integrals can be derived. For example from [23] we
obtain
$$\int_0^{\infty}e^{-x}I_{\nu}(x\sin\;\theta)[\sin(\frac{3\pi\nu}{2})
ber_{2\nu}(2\cos\;\theta\sqrt{ux})+\cos(\frac{3\pi\nu}{2})
bei_{2\nu}(2\cos\;\theta\sqrt{ux})]dx$$
$$=sec\;\theta\sin\;u\;J_{\nu}(u\sin\;\theta),$$
and in particular
$$\int_0^{\infty}e^{-x}I_{2n}(x\sin\;\theta)bei_{4n}(2\cos\;\theta\sqrt{
ux})dx=(-1)^nsec\;\theta\;\sin\;u\;J_{2n}(u\sin\;\theta)$$
$$\int_0^{\infty}e^{-x}I_{2n+1}(x\sin\;\theta)ber_{4n+2}(2\cos\;\theta
\sqrt{ux})dx=(-1)^{n+1}sec\;\theta\;\sin\;u\;J_{2n+1}(u\sin\;\theta).$$
We conclude this section by deriving a further integral representation
for the product $J_{\nu}(ax)J_{\nu}(bx)$ in terms of a $_0F_3$, but
different from (2.4).

Let us consider eq.(2.1) with $\mu=\nu$. According to the quadratic
transformation [10]
$$_2F_1(-m,-\nu-m;\nu+1;\frac{b^2}{a^2})=(1+\frac{b^2}{a^2})^m
\;_2F_1(-\frac{m}{2},\frac{1-m}{2};\nu+1;(\frac{2ab}{a^2+b^2})^2)=$$
$$=(1+\frac{b^2}{a^2})^mm!\sum_{r=0}^{[\frac{m}{2}]}\frac{1}{(m-2r)!r!(\nu
+1)_r}(\frac{ab}{a^2+b^2})^{2r} \eqno(2.27)$$
and (2.1) becomes
$$J_{\nu}(ax)J_{\nu}(bx)=(\frac{1}{4}abx^2)^{\nu}\sum_{m=0}^{\infty}
\frac{(-1)^m(\frac{1}{2}x\sqrt{a^2+b^2})^{2m}}{\Gamma(\nu+m+1)}\cdot$$
$$\cdot\sum_{r=0}^{[\frac{m}{2}]}\frac{(\frac{ab}{a^2+b^2})^{2r}}{(m-2r)!
r!\Gamma(\nu+r+1)}.\eqno(2.28)$$
By using [11]
$$\sum_{m=0}^{\infty}\sum_{r=0}^{[\frac{m}{2}]}c(m,r)=\sum_{m,r=0}^{\infty}
c(m+2r,r)\eqno(2.29)$$
it follows that
$$J_{\nu}(ax)J_{\nu}(bx)=$$ $$\sum_{r=0}^{\infty}\frac{1}{r!\Gamma(\nu+r+1)}(
\frac{abx}{2\sqrt{a^2+b^2}})^{\nu+2r}\sum_{m=0}^{\infty}\frac{(-1)^m(
\frac{1}{2}x\sqrt{a^2+b^2})^{\nu+2r+2m}}{m!\Gamma(\nu+2r+m+1)}=$$
$$=(\frac{1}{2}\frac{abx}{\sqrt{a^2+b^2}})^{\nu}\sum_{r=0}^{\infty}
\frac{(\frac{abx}{2\sqrt{a^2+b^2}})^{2r}}{r!\Gamma(\nu+r+1)}J_{
\nu+2r}(x\sqrt{a^2+b^2})\eqno(2.30)$$
We note, in passing, that (2.30) provides a quick derivation of Weber's
second integral [12]. Indeed, by using
$$\int_0^{\infty}e^{-px^2}x^{\nu+2r+1}J_{\nu+2r}(x\sqrt{a^2+b^2})dx
=\frac{(a^2+b^2)^{\nu/2+r}}{(2p)^{\nu+2r+1}}e^{-\frac{a^2+b^2}{4p}}\eqno(2.31)$$
we have
$$\int_0^{\infty}xe^{-px^2}J_{\nu}(ax)J_{\nu}(bx)dx=\frac{1}{2p}
e^{-\frac{a^2+b^2}{4p}}\sum_{r=0}^{\infty}\frac{(ab/4p)^{\nu+2r}}{r!\Gamma(\nu+r+1)}=$$

$$=\frac{1}{2p}e^{-\frac{a^2+b^2}{4p}}I_{\nu}(\frac{ab}{2p})\eqno(2.32)$$
We also observe that (2.30) can be obtained from the formula [13]
$$_2F_1(a,b;c;x)_2F_1(a,b;c;y)=$$
$$=\sum_{r=0}^{\infty}\frac{(a)_r(b)_r(c-a)_r(c-b)_r}{r!(c)_r(c)_{2r}}(xy)^r
\;_2F_1(a+r,b+r;c+2r;x+y-xy)\eqno(2.33)$$
by applying the confluence principle [14] twice. We first replace x by x/b,
y by y/b, and let b$\rightarrow\infty$. This gives
$$_1F_1(a;c;x)_1F_1(a;c;y)=\sum_{r=0}^{\infty}
\frac{(a)_r(c-a)_r}{r!(c)_r(c)_{2r}}
(-xy)^r\;_1F_1(a+r;c+2r;x+y)\eqno(2.34).$$
Next, we replace x by x/a, y by y/a and let a$\rightarrow\infty$; then
$$\;_0F_1(c;x)\;_0F_1(c;x)=\sum_{r=0}^{\infty}\frac{(xy)^r}{r!(c)_r
(c)_{2r}}\;_0F_1(c+2r;x+y) \eqno(2.35)$$
which is equivalent to the desired result.
If on the right hand side of (2.30) we write [9]
$$J_{\nu+2r}(x\sqrt{a^2+b^2})=$$ $$\frac{2}{\sqrt{\pi}\Gamma(\nu+\frac{1}{2})}
\frac{((1/2)x\sqrt{a^2+b^2})^{\nu+2r}}{(\frac{\nu}{2}+\frac{1}{4})_r
(\frac{\nu}{2}+\frac{3}{4})_r2^{2r}}\int_0^1(1-t^2)^{\nu-\frac{1}{2}
+2r}\cos(xt\sqrt{a^2+b^2})dt,\eqno(2.36)$$
we get ($Re\;\nu>-\frac{1}{2}$)
$$J_{\nu}(ax)J_{\nu}(bx)=\frac{2}{\pi\Gamma(2\nu+1)}(abx^2)^{\nu}
\int_0^1(1-t^2)^{\nu-\frac{1}{2}}\cos(xt\sqrt{a^2+b^2})\cdot$$
$$\cdot\;_0F_3(\nu+1,\frac{\nu}{2}+\frac{1}{4},\frac{\nu}{2}+\frac{3}{4};
\frac{1}{64}a^2b^2x^4(1-t^2)^2)dt\eqno(2.37)$$
and with $b=1$, $x=a^{-1/2}$, and $a=\frac{1}{4}(\sqrt{u^2+2}+\sqrt{u^2-2})^2$,
$$J_{\nu}(\frac{\sqrt{u^2+2}+\sqrt{u^2-2}}{2})J_{\nu}(\frac{
\sqrt{u^2+2}-\sqrt{u^2-2}}{2})=$$ $$\frac{2}{\pi\Gamma(2\nu+1)}\int_0
^1(1-t^2)^{\nu-1/2}\cos(ut)\cdot$$
$$\cdot\;_0F_3(\nu+1,\frac{\nu}{2}+\frac{1}{4},\frac{\nu}{2}+\frac{3}{4}
;\frac{1}{64}(1-t^2)^2)dt \eqno(2.38)$$
Finally, with $\nu=1/2$ (2.37) yields [23]
$$\int_0^1\cos(ut\sqrt{\frac{a^2+b^2}{2ab}})[I_1(u\sqrt{1-t^2})
+J_1(u\sqrt{1-t^2})]\frac{dt}{\sqrt{1-t^2}}$$
$$=(2/u)\sin(u\sqrt{a/2b})\sin(u\sqrt{b/2a}).$$
Further, more complex, evaluations are possible by the same procedure.
\vskip .2in

{\bf 3. The integral $\int_0^{\infty}e^{-\alpha x}J_0(\beta_1
\sqrt{x})J_0(\beta_2\sqrt{x})J_0(\beta_3\sqrt{x})dx$.}

The integral in this section heading, which we denote $I(\beta_1,
\beta_2,\beta_3)$ is an extension of the $\nu=0$ case of Weber's second
integral, eq. (2.32), to which it reduces when one of the three
parameters vanishes.

By writing (see eq.(2.1))
$$J_0(\beta_2\sqrt{x})J_0(\beta_3\sqrt{x})=\sum_{n=0}^{\infty}
\frac{(-1)^n(\frac{1}{2}\beta_3\sqrt{x})^{2n}}{n!^2}\;_2F_1(-n,-n;1;
\frac{\beta_2^2}{\beta_3^2})\eqno(3.1)$$
and [15]
$$\int_0^{\infty}e^{-\alpha
x}x^nJ_0(\beta_1\sqrt{x})dx=\frac{n!}{\alpha^{n+1}}e^{-\beta_1^2/4\alpha}
L_n(\frac{\beta_1^2}{4\alpha}),\eqno(3.2)$$
$L_n(x)$ being a Laguerre polynomial, we first have
$$I(\beta_1,\beta_2,\beta_3)=\frac{1}{\alpha}e^{-\frac{\beta_1^2}{4\alpha}}
\sum_{n=0}^{\infty}\frac{(-1)^n}{n!}(\frac{\beta_3^2}{4\alpha})^n\;_2F_1
(-n,-n;1;\frac{\beta_2^2}{\beta_3^2})L_n(\frac{\beta_1^2}{4\alpha})
\eqno(3.3)$$
Now [16]
$$_2F_1(-n,-n;1;\frac{\beta_2^2}{\beta_3^2})=(1-\frac{\beta_2^2}{\beta_3
^2})^n\;_2F_1(-n,n+1;1;\frac{\beta_2^2}{\beta_2^2-\beta_3^2})$$
$$=(1-\frac{\beta_2^2}{\beta_3^2})^nP_n(\frac{\beta_3^2+\beta_2^2}{\beta
_3^2-\beta_2^2})=\frac{1}{\pi\beta_3^{2n}}\int_0^{\pi}(\beta_2^2+\beta_3^2
-2\beta_2\beta_3\cos\;\theta)^nd\theta.\eqno(3.4)$$
Therefore, in terms of a well known generating function[17], (3.3) becomes
$$I(\beta_1,\beta_2,\beta_3)=\frac{e^{-\frac{\beta_1^2}{4\alpha}}}{\pi\alpha}
\int_0^{\pi}\sum_{n=0}^{\infty}\frac{(1/4\alpha)^n}{n!}
(\beta_2^2+\beta_3^2-2\beta_2\beta_3\cos\;\theta)^n
L_n(\frac{\beta_1^2}{4\alpha})d\theta$$
$$=\frac{e^{-\beta_1^2/4\alpha}}{\pi\alpha}
\int_0^{\pi}e^{-R^2/4\alpha}I_0(\frac{\beta_1R}{2\alpha})d\theta
\eqno(3.5)$$
where $R=(\beta_2^2+\beta_3^2-2\beta_2\beta_3\cos\;\theta)^{1/2}$.
Next, by using Graf's addition theorem [18]
$$I_0(\frac{\beta_1R}{2\alpha})=\sum_{n=0}^{\infty}(2-\delta_{no})
I_n(\frac{\beta_1\beta_2}{2\alpha})I_n(\frac{\beta_1\beta_2}{2\alpha})
\cos(n\theta) \eqno(3.6)$$
and the familiar integral representation
$$\int_0^{\pi}e^{\frac{\beta_1\beta_2}{2\alpha}\cos\;\theta}\cos(n\theta)
d\theta=\pi I_n(\frac{\beta_1\beta_2}{2\alpha}), \eqno(3.7)$$
we finally have
$$\int_0^{\infty}e^{-\alpha
x}J_0(\beta_1\sqrt{x})J_0(\beta_2\sqrt{x})J_0(\beta_3\sqrt{x})dx$$
$$=\frac{1}{\alpha}e^{-\frac{1}{4\alpha}(\beta_1^2+\beta_2^2+\beta_3^2)}
\sum_{n=0}^{\infty}(2-\delta_{0n})F_n(\beta_1,\beta_2,\beta_3)\eqno(3.8)
$$
where
$$F_n(\beta_1,\beta_2,\beta_3)=I_n(\frac{\beta_1\beta_2}{2\alpha})
I_n(\frac{\beta_1\beta_3}{2\alpha})I_n(\frac{\beta_2\beta_3}{2\alpha})
\eqno(3.9)$$

We conclude by sketching the evaluation of the more general integral
$$I_m(\beta_1,\beta_2,\beta_3)=\int_0^{\infty}e^{-\alpha x}J_0(\beta_1
\sqrt{x})J_m(\beta_2\sqrt{x})J_m(\beta_3\sqrt{x})dx \eqno(3.10)$$
where m is a positive integer. In the first place, (3.4) takes the form
$$I_m(\beta_1,\beta_2,\beta_3)=\frac{1}{m!\alpha}(\frac{\beta_2\beta_3}{4\alpha})
^me^{-\frac{\beta_1^2}{4\alpha}}\cdot$$ $$\sum_{n=0}^{\infty}\frac{(-1)^n}{n!}
(\frac{\beta_3^2}{4\alpha})^n\;_2F_1(-n,
-n-m;m+1;\frac{\beta_2^2}{\beta_3^2})L_{m+n}(\frac{\beta_1^2}{4\alpha})
\eqno(3.11)$$
By observing that [19]
$$_2F_1(-n,-n-m;m+1;\frac{\beta_2^2}{\beta_3^2})=(1-\frac{\beta_2^2}{\beta_3^2})^n
\frac{n!(2m)!}{(2m+n)!}C_n^{m+1/2}(\frac{\beta_3^2+\beta_2^2}{\beta_3^2-\beta_2^2})$$
$$=\frac{2^{2m}m!^2}{(2m)!\pi\beta_3^{2n}}\int_0^{\pi}\sin^{2m}\theta\;
(\beta_2^2+\beta_3^2-2\beta_2\beta_3\cos\;\theta)^nd\theta
\eqno(3.12)$$
eq. (3.11) becomes
$$I_m(\beta_1,\beta_2,\beta_3)=\frac{1}{\pi\alpha}\frac{m!}{(2m)!}
(\frac{\beta_2\beta_3}{\alpha})^me^{-\beta_1^2/4\alpha}
\int_0^{\pi}\sin^{2m}\theta\;\sum_{n=0}^{\infty}\frac{(-1/4\alpha)^n}{
n!}\cdot$$
$$(\beta_2^2+\beta_3^2-2\beta_2\beta_3\cos\;\theta)^nL_{m+n}
(\frac{\beta_1^2}{4\alpha})d\theta.\eqno(3.13)$$
Now, from [20],
$$L_{m+n}(x)=\frac{n!}{(m+n)!}e^x(\frac{d}{dx})^me^{-x}x^mL_n^m(x)
\eqno(3.14)$$
and [17]
$$\sum_{n=0}^{\infty}\frac{t^n}{(m+n)!}L_n^m(x)=(xt)^{-m/2}e^t
J_m(2\sqrt{xt}) \eqno(3.15)$$
it follows that (where R has the same meaning as before)
$$\sum_{n=0}^{\infty}\frac{1}{n!}(-\frac{R^2}{4\alpha})^nL_{m+n}
(\frac{\beta_1^2}{4\alpha})=e^{(\beta_1^2-R^2)/4\alpha}2^m(\frac{d}{dx})
^me^{-x}x^m\cdot$$
$$(R\sqrt{\frac{x}{\alpha}})^{-m}I_m(R\sqrt{\frac{x}{\alpha}})|_{x=\beta_1^2/4\alpha}.\eqno(3.16)$$
By inserting (3.16) into (3.13) and using [21]
$$(R\sqrt{\frac{x}{\alpha}})^{-m}I_m(R\sqrt{\frac{x}{\alpha}})=
(\frac{\beta_2\beta_3x}{2\alpha})^{-m}(m-1)!\sum_{n=0}^{\infty}
(m+n)C_n^m(\cos\;\theta)\cdot$$
$$I_{m+n}(\beta_2\sqrt{\frac{x}{\alpha}})I_{m+n}(\beta_3\sqrt{\frac{x}{\alpha}})\eqno(3.18)$$
and [22]
$$\int_0^{\pi}\sin^{2m}\theta\;e^{\frac{\beta_2\beta_3}{2\alpha}\cos\;\theta}C_n^m(\cos\;\theta)d\theta$$
$$=\frac{\pi2^{1-m}(2m+n-1)!}{n!(m-1)!}(\frac{2\alpha}{\beta_2\beta_3})^m
I_{m+n}(\frac{\beta_2\beta_3}{2\alpha}),\eqno(3.18)$$
we finally obtain
$$\int_0^{\infty}e^{-\alpha
x}J_0(\beta_1\sqrt{x})J_m(\beta_2\sqrt{x})J_m(\beta_3\sqrt{x})dx$$
$$=\frac{1}{\alpha}\frac{2^{2m+1}m!}{(2m)!}(\frac{\alpha}{\beta_2\beta_3})^m
e^{-\frac{\beta_2^2+\beta_3^2}{4\alpha}}(\frac{d}{dx})^me^{-x}\cdot
\eqno(3.19)$$
$$\sum_{n=0}^{\infty}\frac{(n+m)(2m+n-1)!}{n!}I_{m+n}(\beta_2\sqrt{\frac{x}{\alpha}})I_{m+n}(\beta_3\sqrt{\frac{x}{\alpha}})I_{m+n}(\frac{\beta_2\beta_3}{2\alpha})|_{x=\beta_1^2/4\alpha}.$$
In particular, by letting $\beta_3\rightarrow0$, after dividing by
$\beta_3^m$,
$$\int_0^{\infty}e^{-\alpha
x}J_0(\beta_1\sqrt{x})J_m(\beta_2\sqrt{x})x^{m/2}dx$$
$$=\frac{1}{\alpha}e^{-\beta_2^2/4\alpha}(\frac{d}{dx})^m
e^{-x}(\frac{x}{\alpha})^{m/2}I_m(\beta_2\sqrt{\frac{x}{\alpha}})|
_{x=\beta_1^2/4\alpha}\eqno(3.20)$$
$$=\frac{1}{\alpha}(\frac{\beta_2}{2\alpha})^me^{-\frac{\beta_1^2+\beta_2^2}{4\alpha}}
\sum_{n=0}^m(-1)^n{m\choose n}(\frac{\beta_1}{\beta_2})^nI_n(
\frac{\beta_1\beta_2}{2\alpha}).$$
which is of interest in connection with formulas (39)-(42) on page 186
of reference [15].

In conclusion, we point out that the derivation of (2.24) can be
extended to give

$$\int_0^{\infty}e^{-\beta x}\;_0F_3(\mu , \nu , \nu+\frac{1}{2};
-a^2x^2)dx=(2a)^{1-\mu}\Gamma(\mu)\Gamma(2\nu)\beta^{\mu-2\nu-1}
J_{\mu-1}(4a/\beta) \eqno(3.21)$$
Therefore, the results of section 2, for example, are capable of
extension in a variety of directions. We leave this for the future and
merely quote one example, interesting because it contains each of the
four types of Bessel functions:
$$\int_0^{\infty} xJ_1(ax)I_1(ax)Y_0(x)K_0(x)\; dx=-(2\pi a^2)^{-1}
\ln(1-a^4)\eqno(3.22)$$
where $0<a<1$.

\centerline{\bf References}
\vskip .2in
\noindent
[1] T.J. Lardner, Siam Review \underbar{11}, 69 (1969).

\noindent
[2] Bateman Manuscript Project, Higher Transcendental Functions
(McGraw-Hill, 1953), Vol. 2, p.11, eq. (47).

\noindent
[3] Ibid. Vol. 1, p.64, eq. (23).

\noindent
[4] Ibid. Vol. 2, p.52, eq. (31).

\noindent
[5] Ibid. Vol. 2, p.51, eq. (29).

\noindent
[6] L.J. Slater, Generalized Hypergeometric Functions (Cambridge 1966),
Sec. 4.8.

\noindent
[7] I.N. Sneddon, The Use of Integral Transforms (McGraw-Hill 1972),
Chap. 5, p.309.

\noindent
[8] See Ref.[2] Vol. 2, p.96, eq. (59).

\noindent
[9] Ibid. Vol. 2, p.81, eq. (10).

\noindent
[10] See Ref. [1], p. 110, eq. (1).

\noindent
[11] H.M. Srivastava and H.L. Manocha, A treatise on Generating
Functions (John Wiley 1983), p. 101, eq. (6).

\noindent
[12] See Ref. [2], Vol. 2, p.50, eq. (25).

\noindent
[13] J.L. Burchnall and T.W. Chaundy, Quart. J. Math. Oxford Ser. 11,
249-270 (1940).

\noindent
[14] See Ref. [11], p.36.

\noindent
[15] Bateman Manuscript Project, Tables of Integral Transforms
(McGraw-Hill 1953), Vol. 2, p. 30, eq. (13).

\noindent
[16] See Ref. [1], p. 157, eq. (15).

\noindent
[17] See Ref. [2], p. 189, eq. (18).

\noindent
[18] Ibid. Vol. 2, p.44, eq.(5).

\noindent
[19] Ibid. Vol. 2, p.177, eq.(31).

\noindent
[20] Ibid. Vol. 2, p.190, eq.(28).

\noindent
[21] Ibid. Vol. 2, p.101, eq.(30).

\noindent
[22] See Ref. [15], p.281, eq.(7).

\noindent
[23] A.P. Prudnikov, Ya. A. Brychkov and O.I. Marichev, Integrals and
Series (Gordon and Breach, N.Y. 1989) Vol. 3 Sec. 7.16.3 Eq. (6).

\noindent
[24] Ibid. Sec. 7.16.2 Eq.(15).

\end{document}